\documentstyle[11pt,titlepage,amssymb,euscript]{article}
\newtheorem{theorem}{Theorem}[section]
\newtheorem{proposition}[theorem]{Proposition}
\newtheorem{lemma}[theorem]{Lemma}
\newtheorem{corollary}[theorem]{Corollary}
\newtheorem{definition}[theorem]{Definition}
\newtheorem{example}[theorem]{Example}
\newtheorem{construction}[theorem]{Construction}

\newcommand{\R}{\mbox{${\Bbb R}$}}

\newcommand{\Z}{\mbox{${\Bbb Z}$}}

\newcommand{\D}{\mbox{$\nabla$}}

\newcommand{\g}{\mbox{$\gamma$}}

\newcommand{\ns}{\normalshape} 
\newcommand{\epi}{\rightarrow\hspace{-9pt}\rightarrow}
\newcommand{\p}{\mbox{$\pi$}}
\newcommand{\scr}{\EuScript}
\newcommand{\th}{\mbox{$\tilde h$}}
\newcommand{\tp}{\mbox{$\tilde\pi$}}
\newcommand{\ve}{\mbox{$\varepsilon$}}
\newcommand{\vp}{\mbox{$\varphi$}}

\newcommand{\wD}{\mbox{$\widetilde{\nabla}$}}
\newcommand{\wM}{\mbox{$\widetilde{M}$}}

\renewcommand{\H}{\mbox{${\Bbb H}$}}
\renewcommand{\P}{\mbox{${\Bbb P}\hspace*{1pt}$}}

\renewcommand{\b}{\mbox{$\beta$}}
\renewcommand{\l}{\mbox{$\lambda$}}

\makeatletter
\@addtoreset{equation}{section}
\makeatother

\font\heads = cmbx12

\title{\Large\bf SPACES OF GEODESICS:\\PRODUCTS, COVERINGS, CONNECTEDNESS}

\author{
John K. Beem\\
\noalign{\medskip}
{\small Department of Mathematics}\\
{\small University of Missouri}\\
{\small Columbia MO \ 65211}\\
{\small USA}\\
{\small mathjkb@mizzou1.missouri.edu}
\and Robert J. Low\\
\noalign{\medskip}
{\small Department of Mathematics}\\
{\small Coventry University}\\
{\small Coventry \ CV1 5FB}\\
{\small UK}\\
{\small roblow@coventry.ac.uk}
\and Phillip E. Parker\\
\noalign{\medskip}
{\small Department of Mathematics}\\
{\small The Wichita State University}\\
{\small Wichita KS \ 67260-0033}\\
{\small USA}\\
{\small pparker@twsuvm.uc.twsu.edu}
\vspace*{2cm}
}
\date{ January 14, 1995 \\
\vspace*{2cm}
}

\begin{document}

\maketitle

\begin{abstract}
We continue our study of the space of geodesics of a
manifold with linear connection. We obtain sufficient conditions for a product
to have a space of geodesics which is a manifold. We investigate the
relationship of the space of geodesics of a covering manifold to that of the
base space. We obtain sufficient conditions for a space to be
geodesically connected in terms of the topology of its space of geodesics.
We find the space of geodesics of an $n$-dimensional Hadamard manifold is
the same as that of $\R^n$.
\end{abstract}

\section{\heads Introduction}

A topological space is {\it locally Euclidean\/} if and only if
every point has a neighborhood homeomorphic to an open set in some $\R^n$.
Such a space is sometimes called a manifold which is not Hausdorff.  Here,
we shall assume that locally Euclidean spaces are paracompact (not
necessarily Hausdorff) and connected, and a {\it manifold\/} will be a
Hausdorff locally Euclidean space.  Note that the theory of differential
structures (smoothness) does not depend on the Hausdorff condition, so that
one may speak of smooth ($C^\infty$) locally Euclidean spaces.  We agree
that all locally Euclidean spaces are smooth.  Thus there is a tangent
bundle $TE$, etc.  Notations will be given explicitly only for manifolds,
but may be used for locally Euclidean spaces whenever appropriate.

For a smooth manifold $M$ with tangent bundle $TM$, let $T'\hspace*{-1pt}M$
denote the reduced tangent bundle consisting of $TM$ less the zero section
and let $PM$ denote the quotient space of $T'\hspace*{-1pt}M$ obtained by
identifying proportional vectors at each $p\in M$.  If $M$ is
$n$-dimensional, then $PM$ has compact fibers diffeomorphic to the real
projective space $\P^{n-1}$.  We shall regard points in $PM$ as tangent
lines to curves in $M$.  There is also the double covering $SM$ of $PM$ in
which the proportionality must be by positive numbers and the compact
fibers are diffeomorphic to the sphere $S^{n-1}$.

When $M$ is provided with a linear connection, there are geodesics in $M$.
We agree that all geodesics are to be inextendible (maximally extended).
Define the {\it space of (unoriented) geodesics\/} $G(M)$ to be the
quotient of $PM$ obtained by identifying points that lie on the lift of a
common geodesic.  Elements of $G(M)$ may be regarded as equivalence classes
$[\g]$ of parameterized geodesics $\g :(a,b)\longrightarrow M$.  Similarly,
there is the space of {\it oriented\/} geodesics $G^+(M)$, the
quotient of $SM$.  Observe that different connections on $M$ may yield
topologically different spaces of geodesics.

One usually regards {\em closed\/} and {\em periodic\/} as synonymous for
geodesics.  However, it is possible for a null geodesic to retrace its path
without the tangent vector returning to itself:  since it has zero length,
it need return only to a scalar multiple of itself.  A linear connection
need not be metric, and then this could happen with all the geodesics which
retrace themselves.  Thus we shall distinguish between {\em closed\/}
geodesics in which the image is a closed curve and {\em periodic\/}
geodesics in which the image of the natural lift to $TM$ is a closed curve.
Observe that for spacelike and timelike geodesics of a pseudoriemannian
metric tensor, closed and periodic are still synonymous.  In particular,
this is true for all geodesics of a Riemannian metric.

This paper continues \cite{BP2} and \cite{L1}, and contains proofs of most
of the results announced in \cite{Pr}; proofs of the others will appear later.

Let $(M_1,\D_1)$ and $(M_2,\D_2)$ be two manifolds with linear connections.
In Section 2 we consider the space of geodesics of the
product $(M_1\times M_2, \D_1\times \D_2)$.  We find that if $G(M_1\times
M_2)$ is locally Euclidean, then both $G(M_1)$ and $G(M_2)$ must be
locally Euclidean.  On the other hand, if both $G(M_1)$
and $G(M_2)$ are locally Euclidean, then $G(M_1\times M_2)$ is locally
Euclidean if and only if both $(M_1, \D_1)$ and $(M_2, \D_2)$ fail to have
any closed geodesics.  In particular, one may easily
construct examples with both $G(M_1)$ and $G(M_2)$ manifolds,
but $G(M_1\times M_2)$ not a manifold.  Stronger results may be
obtained when the geodesics of
$(M_1\times M_2, \D_1\times \D_2)$ satisfy the nonreturning property.
Under this assumption, $G(M_1\times M_2)$ is a manifold if and only if both
$G(M_1)$ and $G(M_2)$ are manifolds.

In Section 3 we consider smooth coverings with
projections that preserve geodesics.  If $p: \wM \epi M$ is
such a covering, then there is an induced
map $\hat{p} : G(\wM) \rightarrow G(M)$ and this induced map
commutes with the natural projections $T'(M)\epi G(M)$ and $T'(\wM) \epi
G(\wM)$. In general $G(\wM)$ may
fail to be a manifold when $G(M)$ is a manifold.
However, if $(M, \D)$ is nonreturning
then $G(\wM)$ will be a manifold when $G(M)$ is a manifold.

The Hopf-Rinow Theorem guarantees that if $(M, g)$ is a
complete Riemannian (positive definite) manifold, then
each pair of points may be joined by at least one geodesic.
Thus, one says that a complete
Riemannian manifold is geodesically connected.  In spaces
which are not positive definite, there is no
Hopf-Rinow Theorem and the problem of deciding if
each pair of points may be joined by some geodesic
can be quite difficult.  In Section 4 we consider
the problem of geodesic connectedness.  We find
that $(M, \D)$ is geodesically connected if it is
unitrace and $G(M)$ is
Hausdorff.   We define sky bundles and find that
if one has a sky bundle with the top Stiefel-Whitney
class nontrivial, then one may connect pairs of
points with geodesic segments.  This yields a
new approach to the problem of deciding when a
manifold with a linear connection is geodesically connected.
In Section 4, we also show that for Hadamard manifolds
the geodesic space is the same as
that of $\R^n$.

\section{\heads Products}

We begin by recalling some definitions and results from \cite{BP1,BP2} and
\cite{Pr}. From \cite{BP2} we need
\begin{definition}
$(M,\D)$ is {\em unitrace\/}
if and only if for each $p\in M$ and each neighborhood $V$ of $p$ there is
a simple convex neighborhood $U$ of $p$ with $U\subseteq V$ such that any
geodesic which enters $U$ either leaves and never returns or retraces the
same path every time it does return.  If each point has arbitrarily small
neighborhoods such that no geodesic ever returns, we say $(M,\D)$ is {\em
nonreturning.}
\end{definition}
Observe that nonreturning implies unitrace
implies the geodesic foliation of $PM$ is regular in the sense of Palais
\cite{P}. As in \cite{Pr}, we make the
\begin{definition}
$(M,\D)$ is said to be {\em geodesically regular} if and only if the geodesic
foliation of $PM$ is regular.
\label{gr}
\end{definition}
This means that each point in $PM$ has a neighborhood which is intersected by
each lifted geodesic at most once. In particular, this implies that each
point of $M$ has a neighborhood which every geodesic image
crosses at most finitely many times. This is not true on the Misner cylinder
$\scr M$
\cite[p.\,177ff\,]{HE}. When the connection is flat, the crossings cannot
be parallel and arbitrarily close. This is not true on the flat M\"obius band
$\frak M$. From \cite[p.\,19]{P} we obtain \cite[Theorem 3.15]{Pr}
\begin{theorem} $(M,\D)$ is geodesically regular if and only if $G(M)$ is
locally Euclidean.  $\Box$
\label{grle}
\end{theorem}

Finally we recall \cite{BP1,BP2} \begin{definition} We say that $(M,\D)$ is
{\em disprisoning\/} if and only if each end of each geodesic eventually
leaves each compact set.  $(M,\D)$ is {\em $\D$-pseudoconvex\/} if and only
if for every compact set $K$ there is a compact set $K'$ such that every
geodesic segment with both endpoints in $K$ lies entirely in $K'$.
\end{definition}
One regards pseudoconvexity as a sort of interior completeness condition.

\medskip
Now let $(M_1,\D_1)$ and $(M_2,\D_2)$ be manifolds with linear connections.
In general, the space of geodesics of the product cannot readily
be determined from those of the factors.  Some properties, however, are
reasonably well shared by a product and its factors.
\begin{proposition}
Let $(M_1,\D_1)$ and $(M_2,\D_2)$ be geodesically regular.  The product
$(M_1\times M_2,\D_1\times\D_2)$ is geodesically regular if and only if
neither factor has a closed geodesic.
\label{grncd}
\end{proposition}
{\bf Proof:} Assume one factor, say $M_1$, has a closed geodesic \g\ with
$\g'(0) = c\g'(a)$.  Let \b\ be any geodesic of $M_2$ with 0 in the domain
of \b.  Set $\l(\ve,t) = (\g(t),\b(\ve t))$.  For each fixed \ve, one has
$\l(\ve,t)$ a geodesic of $M_1\times M_2$ and $\l(\ve,a) \rightarrow
(\g(0),\b(0)) = (\g(a),\b(0))$ as $\ve\rightarrow 0$.  Furthermore,
$\l'(\ve,0)$ and $\l'(\ve,a)$ converge to $(\g'(0),0)$ and $(\g'(a),0)$,
respectively, as $\ve\rightarrow 0$.  It follows that each neighborhood of
the point in $P(M_1\times M_2)$ corresponding to $(\g'(0),0)$ is
intersected more than once by the curve of tangent lines corresponding to
$\l'(\ve,t)$ for all sufficiently small $\ve\neq 0$.  Thus the product is
not geodesically regular.

Conversely, if $M_1\times M_2$ is not geodesically regular, then there
exists $v = (v_1,v_2)\in T_x(M_1\times M_2)\cong T_{x_1}M_1\times
T_{x_2}M_2$ where geodesic regularity fails.  Let $U$ be any simple convex
normal neighborhood of $x$.  Then there is a sequence of geodesics $\l_n =
(\g_n,\b_n)$ such that $\l'_n(0)\rightarrow v$ and the tangent lines
containing $\l'_n(t_n)$ converge to the tangent line containing $v$ for
some sequence $\{ t_n\}$ where $\l_n$ leaves $U$ between 0 and $t_n$.  Let
$U_i$ be the projection of $U$ onto $M_i$ for $i = 1,2$.  Without loss of
generality, we may as well assume that the geodesics $\g_n$ of $M_1$ leave
$U_1$ after $t=0$ and return before $t=t_n$ with tangent lines converging
to that containing $v_1$.  Since $M_1$ is geodesically regular, it follows
that $\g_n$ must retrace its image in $U_1$ for all large $n$.  Therefore
$M_1$ has a closed geodesic.  ${}_\Box$

\medskip\noindent
Using Theorem~\ref{grle}, we obtain
\begin{theorem}
Let $G(M_1)$ and $G(M_2)$ be locally Euclidean.  The space of geodesics of
the product $G(M_1\times M_2)$ is locally Euclidean if and only if neither
factor has a closed geodesic.  ${}_\Box$
\end{theorem}

The space of geodesics of the product need not be Hausdorff either, even if
that of each factor is.  A simple example is the standard pseudoeuclidean
$\R^n\times S^k$, $n,k>0$.

\begin{proposition}
$(M_1\times M_2,\D_1\times\D_2)$ is nonreturning if and
only if $(M_1\times M_2,\D_1\times\D_2)$ is unitrace if and only if
$(M_1,\D_1)$ and $(M_2,\D_2)$ are both nonreturning.
\label{nub}
\end{proposition}
{\bf Proof:} Since geodesics of the product are ordered pairs of geodesics
of the factors, it follows easily that the product is nonreturning if and
only if both factors are.

Clearly, if the product is nonreturning then it is unitrace.  It only
remains to show that if the product is unitrace, then it is nonreturning.
To this end, assume that the product has a geodesic $\g = (\g_1,\g_2)$ that
leaves a neighborhood $U$ of $\g(0)$ and later returns and retraces its
image.  At least one of $\g_1$ or $\g_2$ is nontrivial, so we may as well
assume that $\g_1$ is nontrival.  Then there is some $a$ and $c$ with
$\g'_1(0) = c\g'_1(a)$.  Let \b\ be any nontrivial geodesic of $(M_2,\D_2)$
and set $\l(\ve,t) = (\g_1(t),\b(\ve t))$.  Letting $\ve\rightarrow 0$ as
in the first part of the proof of Proposition~\ref{grncd}, we obtain a
failure of geodesic regularity.  But the product is unitrace, hence
geodesically regular; contradiction.  ${}_\Box$

\begin{theorem}
If $G(M_1\times M_2)$ is locally Euclidean, then so are
$G(M_1)$ and $G(M_2)$.
\end{theorem}
{\bf Proof:} Using
Theorem~\ref{grle}, we need only show that geodesic regularity of the
product implies geodesic regularity of the factors.  Fixing $(x_1,x_2) \in
M_1\times M_2$, the result follows easily from the identification of $M_1$
with the submanifold $M_1\times\{ x_2\}$ of $M_1\times M_2$ and of $M_2$
with the corresponding $\{ x_1\}\times M_2$.  ${}_\Box$

\begin{lemma}
$(M_1\times M_2, \D_1\times\D_2)$ is pseudoconvex if and only
if $(M_1,\D_1)$ and $(M_2,\D_2)$ both are.
\end{lemma}
{\bf Proof:} Let
$r_1$ and $r_2$ be the canonical projections onto the first and second
factors, respectively, and let $(x_1,x_2)$ be a fixed point of the product.

First, assume that the product is pseudoconvex and let $K$ be a compact set
in $M_1$.  Then $K\times \{ x_2\}$ is a compact subset of the product,
hence there is a compact $H \subseteq M_1\times M_2$ such that any geodesic
segment of $M_1\times M_2$ with endpoints in $K\times\{ x_2\}$ must lie in
$H$.  If $\g:[a,b]\rightarrow M_1$ is a geodesic segment in $M_1$ with
endpoints in $K$, then $\l = (\g,x_2)$ is a geodesic segment in $M_1\times
M_2$ with endpoints in $K\times\{ x_2\}$ and image in $H$.  Consequently,
$\g = r_1\circ\l$ lies in the compact set $r_1(H)$ for any such \l.  It
follows that $(M_1,\D_1)$ is pseudoconvex; the argument for $(M_2,\D_2)$ is
similar.

Conversely, assume that the factors are pseudoconvex.  Now let $K\subseteq
M_1\times M_2$ be compact and set $K_1 = r_1(K)$ and $K_2 = r_2(K)$.  There
are compact $H_i\subseteq M_i$ such that any geodesic segment of $M_i$ with
endpoints in $K_i$ must lie in $H_i$.  Now $H = H_1\times H_2$ is compact.
Since geodesics of the product are ordered pairs of geodesics of the
factors, it follows that any geodesic of $M_1\times M_2$ with endpoints in
$K\subseteq K_1\times K_2$ must lie in $H$.  Therefore the product is
pseudoconvex.  ${}_\Box$

\medskip\noindent
Using this lemma and Corollary 5.6 of \cite{BP2}, we obtain
\begin{theorem}
If $(M_1\times M_2,\D_1\times\D_2)$ is nonreturning,
then these are equivalent:
 \begin{enumerate}
 \item $G(M_1\times M_2)$ is Hausdorff;
 \item $G(M_1\times M_2)$ is a manifold;
 \item $(M_1\times M_2, \D_1\times\D_2)$ is pseudoconvex;
 \item $G(M_1)$ and $G(M_2)$ are Hausdorff;
 \item $G(M_1)$ and $G(M_2)$ are manifolds;
 \item $(M_1,\D_1)$ and $(M_2,\D_2)$ are pseudoconvex.  ${}_\Box$
 \end{enumerate}
\label{mct}
\end{theorem}

Combining Proposition \ref{nub} and Theorem \ref{mct}, we obtain
sufficient
conditions for a product to have a space of geodesics which is a manifold.
\begin{theorem}
If $(M_1,\D_1)$ and $(M_2,\D_2)$ are each both nonreturning and pseudoconvex,
then $G(M_1\times M_2)$ is a manifold. ${}_\Box$
\end{theorem}

The general problem to construct $G(M_1\times M_2)$ from $G(M_1)$ and
$G(M_2)$ is complicated and somewhat difficult, and more work needs to be
done here.  We shall now construct $G(\R \times M)$ for $M$ with a linear
connection.
\begin{construction}
{\normalshape
Coordinatize $R \times M$
by $(s,x)$ where $s \in R$ and $x$ gives the coordinates of a point in $M$.
We shall call $\partial_s$ the horizontal direction and refer to directions
tangent to $M$ as vertical.  Let $p$ be the natural projection to $M$, and
identify $M$ with $\{0\} \times M$.

First, consider all the geodesics in $\R\times M$ that project to the
nondegenerate geodesic in $\R$.  All geodesics of this form must pass
through some point $(0,x)$ for some $x$, and are completely specified by
their tangent vector there.  We normalize so that the projection to $\R$
has unit velocity.  Then the nonvertical geodesics are given exactly by
$TM$.  Next, there are the vertical geodesics to consider.  For each $s \in
\R$, we obtain the whole of $G(M)$.  Thus, as a set, we have $TM \sqcup \R
\times G(M)$.  The only remaining question is, what is the topology?  This
can be settled by seeing what sequences in $TM$ converge to which points
(if any) in $\R \times G(M)$.

Recall that each point $(x,v)$ of $TM$ specifies the geodesic through
$(0,x)$ with tangent $(1,v)$.  Let $(x_n,v_n)$ be a sequence in $TM$; under
what circumstances will this sequence, regarded as a set of geodesics in
$\R \times M$, converge to some vertical geodesic?  Denote the geodesic
corresponding to $(x_n,v_n)$ by $\gamma_n$, where $\gamma(0) = (0,x_n)$,
and $\dot{\gamma}_n = (1,v_n)$, and let $\gamma$ be the vertical geodesic
with $\gamma(0) = (s,x)$ and $\dot{\gamma}_n = (0,v)$.  Then $\gamma_n
\rightarrow \gamma$ if and only if $p(\gamma_n) \rightarrow \gamma$, and
the point of intersection of $\gamma_n$ with $\{s\} \times M$ converges to
$(s,x)$.  ${}_\Box$
}
\end{construction}

\medskip\noindent
This construction can be iterated to describe the set
$G(\R^k\times M)$ for any standard pseudoeuclidean $\R^k$, but the details
of the topology are rather complicated.

In the case of $\R\times M$, it is easy to see that if $M$ has a closed
geodesic, then a sequence in $TM$ converging to one vertical copy of this
geodesic will converge to all copies as for the classical cylinder.  It is
also easy to carry out this construction explicitly for the case where $M$
is $\R$ and obtain a M\"obius band.

\begin{example}
{\normalshape
First, consider all the nonvertical geodesics
in $\R\times \R = \{ (x,y)\}$.  These geodesics are completely described by
$T\R$ {\it via\/} the $y$-intercept $(0,y)$ and the normalized velocity
$(1,v)$.  Now $T\R$ is another copy $\R\times \R = \{ (u,v)\}$.  The point
$(u,v)$ corresponds to the nonvertical geodesic with $y$-intercept $(0,u)$
and normalized velocity $(1,v)$.

Next we consider those sequences of nonvertical geodesics which converge to
a vertical geodesic, in order to attach the boundary copy of \R\ in the
correct way.  Any sequence of nonvertical geodesics through the point
$(x,0)$ in $\R^2_{xy}$ which converges to the vertical geodesic through
$(x,0)$ has a sequence of normalized velocities $(1,v_n)$ with $v_n
\rightarrow \pm\infty$.  This yields the sequence of points $(-xv_n, v_n)$
in $\R^2_{uv}$.  Observe that these points lie along the line through the
origin of slope $-1/x$, interpreted as the $v$-axis when $x = 0$.

Now shrink $\R^2_{uv}$ to the unit disc; {\it e.g.,} map the point with
polar coordinates $(r,\theta)$ to the point $(\tanh r, \theta)$.  We see
that we should attach one boundary point to each antipodal position pair
{\it except\/} the pair $(\pm 1,0)$ at opposite ends of the $u$-axis.  This
recovers $G(\R^2)$ as the projective plane less a point:  the open M\"obius
band. ${}_\Box$
}
\end{example}

\medskip
If we consider instead the oriented geodesics of $\R\times\R$, we
obtain two copies of the disc glued along the boundaries in the oriented
way except at the points $(\pm 1,0)$ again.  This yields $S^2$ less two
antipodal points, diffeomorphic to $TS^1$.

\section{\heads Coverings}

Let $p:  \widetilde{M} \epi M$ be a smooth covering map; in
particular, we assume that $p$ is a local diffeomorphism which evenly
covers $M$.  We have
\begin{lemma}
If $p$ is a such a covering map, then $p$ and the induced tangent map $p_*$
are open. ${}_\Box$
\end{lemma}

\begin{definition}
We say that $p:  (\widetilde{M},\widetilde{\nabla})
\epi (M,\D)$ is a {\em geodesic\/} covering map if and only if
$p$ takes the image of each geodesic in $(\widetilde{M},
\widetilde{\nabla})$ onto the image of a geodesic
in $(M, \nabla)$.
\end{definition}
Such maps preserving geodesics up to parameterization are
called {\em projective} ({\em cf.} \cite{K,KN} for older terminology).

Now suppose $p$ is a geodesic covering map.  Let $\p :  T'M \epi G(M)$, and
$\tp :  T'\wM \epi G(\wM)$ be the natural projections.
\begin{proposition}
There is an induced map $\hat{p} :  G(\wM)\rightarrow
G(M)$ which is a continuous, open surjection (thus a quotient map) such
that $\hat{p}\tp = \p p_*$.
\end{proposition}
{\bf Proof:} Fix $[\g] \in
G(\wM)$ and pull back to $T'\wM$ by using $G = \tp^{-1}([\g])$.  We shall
show that if $u,v\in G$, then $\p p_*(u) = \p p_*(v)$, which will show that
for each $[\g]\in G(\wM)$ we have a uniquely defined element $\hat{p}([\g])
= \p p_*(G) = \p p_*(u) = \p p_*(v) \in G(M)$.  If $u,v\in G$, then there
is a geodesic \b\ in \wM\ with $\b\in[\g]$ and $\b'(0) = u$ and $\b'(a) =
cv$.  Using the fact that $p$ is a geodesic covering, it follows that
$p\circ\b$ has tangent vectors proportional to $p_*u$ and $p_*v$ at 0 and
$a$, respectively.  Consequently, \p\ takes $p_*u$ and $p_*v$ to the same
image in $G(M)$.  Note that $\hat{p}[\g] = [p\circ\g]$, and $\hat{p}$ is
surjective since $p$ is.  By construction, $\hat{p}\tp = \p p_*$.

If $W$ is open in $G(M)$, then $p_*^{-1}(\p^{-1}(W))$ is open in $T'\wM$
since \p\ and $p_*$ are continuous.  Thus $\tp(p_*^{-1}(\p^{-1}(W)))$ is
open in $G(\wM)$ because \tp\ is an open map by Lemma 3.1 of \cite{BP2}.
On the other hand, $\hat{p}^{-1}(W) = \tp(p_*^{-1}(\p^{-1}(W)))$ which
establishes the continuity of $\hat{p}$.

If $U$ is open in $G(\wM)$, then the continuity of \tp\ yields
$\tp^{-1}(U)$ open in $T'\wM$.  Both $p_*$ and \p\ are open, so we find
that $\hat{p}(U) = \p(p_*(\tp^{-1}(U)))$ is open in $G(M)$.  Therefore
$\hat{p}$ is an open map.  ${}_\Box$

\begin{lemma}
If $M$ is disprisoning, unitrace, or nonreturning, then so is
$\widetilde{M}$, respectively.
\label{dunc}
\end{lemma}
{\bf Proof:} If a
geodesic in \wM\ has one end imprisoned in a compact set, so does its
projection in $M$.  Thus if \wM\ is not disprisoning, then neither is $M$.
Similarly, if \wM\ is not unitrace or nonreturning, neither is $M$. ${}_\Box$

\begin{proposition}
If $M$ is pseudoconvex and nonreturning, then so is \wM.
\label{pnc}
\end{proposition}
{\bf Proof:} It suffices to show that
\wM\ is pseudoconvex.  Fix an auxiliary complete Riemannian (positive
definite) metric tensor $h$ on $M$ and lift it to \th\ on \wM.  Observe
that $p$ is also a geodesic covering for \th\ and $h$ by construction.  It
follows that \th\ is complete.  If $\g:(a,b)\rightarrow\wM$ is a geodesic
of $(\wM,\wD)$, then its \th -length is equal to the $h$-length of the
image $p\circ\g$ in $M$.  Assume that \wM\ fails to be \wD -pseudoconvex
because there are geodesic segments with endpoints in the compact set
$K\subseteq\wM$ which fail to be contained in any compact set of \wM.  The
continuity of $p$ yields that the image $p(K)$ is a compact set in $M$.
Since $M$ is \D -pseudoconvex, there is a compact set $H$ of $M$ such that
any \D -geodesic of $M$ with endpoints in $p(K)$ must remain in $H$.  Cover
$H$ with a finite number of \D -convex normal neighborhoods having compact
closure such that any \D -geodesic that leaves one of these neighborhoods
fails to return.  For each such neighborhood, there is a maximal $h$-length
for each \D -geodesic segment lying in that neighborhood.  Clearly, there
is a real number $L$ such that each \D -geodesic segment which lies in $H$
has $h$-length less than $L$.  Now choose any \wD -geodesic segment $\g_1$
with \th -length greater than $L$ and with endpoints in $K$.  (Such a \wD
-geodesic segment exists since \th\ is complete and no compact set contains
all the \wD -geodesic segments with endpoints in $K$.)  Note that
$p\circ\g_1$ is a \D -geodesic segment in $M$ which lies in $H$ and has
$h$-length greater than $L$; but this contradicts the definition
of $L$. ${}_\Box$

\medskip\noindent
Using this together with Corollary 5.6 of \cite{BP2}, we
immediately obtain
\begin{corollary}
If $M$ is nonreturning and $G(M)$ is a
manifold, then $G(\wM)$ is also a manifold.  ${}_\Box$
\end{corollary}

The next example shows that one may not weaken the assumption that $(M,\D)$
is nonreturning to unitrace.  It also shows shows that for $(M,\D)$
unitrace, $G(M)$ Hausdorff is not equivalent to $(M,\D)$ being pseudoconvex.
\begin{example}
{\normalshape
Let $M_0$ be the real projective plane $\P^2$
with the usual elliptic metric and connection.  If $M$ is $M_0$ less a
point, then it is not hard to show that $G(M) = G(M_0) = \P^2$ is a
manifold and $(M,\D)$ is unitrace but not pseudoconvex.  The universal
cover $(\wM,\wD)$ is $S^2$ less two antipodal points with the usual
connection.  Since $M$ is $\P^2$ less a point, one may select a geodesic
\g\ which is not closed.  There are two geodesics $\tilde{\g}_1$ and
$\tilde{\g}_2$ of $(\wM,\wD)$ lying over \g\ and one may construct a sequence
of closed geodesics of $(\wM,\wD)$, the images of which have a Hausdorff
closed limit equal to the union of the images of $\tilde{\g}_1$ and
$\tilde{\g}_2$.  Thus $G(M)$ is Hausdorff but $G(\wM)$ is not. ${}_\Box$
}
\label{unp}
\end{example}

\begin{proposition}
Assume that $p$ is a finite covering.  Then \wM\ is
pseudoconvex if and only if $M$ is.
\end{proposition}
{\bf Proof:} First
assume that $(\wM,\wD)$ is pseudoconvex and let $K$ be a compact set in
$M$.  Since $p$ is a finite cover, $p^{-1}(K)$ is also compact.  Let
$\widetilde{H}$ be a compact set in \wM\ such that any \wD -geodesic
segment with endpoints in $p^{-1}(K)$ lies in $\widetilde{H}$.  Each \D
-geodesic segment in $M$ with endpoints in $K$ must lift to a \wD -geodesic
segment with endpoints in $p^{-1}(K)$, hence the lift must lie in
$\widetilde{H}$.  It follows that $(M,\D)$ is pseudoconvex because each \D
-geodesic with endpoints in $K$ must lie in the compact set
$p(\widetilde{H})$.

Conversely, assume that $(M,\D)$ is pseudoconvex and let $\widetilde{K}$ be
a compact set in \wM.  Since $p(\widetilde{K})$ is compact, there is a
compact set $H$ in $M$ such that any \D -geodesic segment with endpoints in
$p(\widetilde{K})$ lies in $H$.  Because $p$ is a finite covering,
$p^{-1}(H)$ is a compact set.  Now, each \wD -geodesic segment with
endpoints in $\widetilde{K}$ must lie in $p^{-1}(H)$ because it projects to
a \D -geodesic segment with endpoints in $p(\widetilde{K})$ and lying in
$H$.  ${}_\Box$

\section{\heads Connectedness}

It will be convenient to modify the definition of sky as used by Low
\cite{L1,L2}.
\begin{definition}
For any point $x\in M$, the set $X$ of
all geodesics through $x$ is called the {\em sky} of $x$.
\end{definition}
We shall follow this upper and lower case convention to denote skies.

\begin{theorem}
If $(M,\D)$ is unitrace and $G(M)$ is Hausdorff, then $M$
is geodesically connected.
\end{theorem}
{\bf Proof:} First, we observe
that $G(M)$ is a manifold by Proposition 5.1 of \cite{BP2}.  It is easy to
see that the sky of a point $x\in M$ gives a subset $X$ of $G(M)$ that is
diffeomorphic to $\P^{n-1}$.

Because $M$ has a linear connection, $x$ has a convex normal neighborhood.
Furthermore, since $M$ is unitrace, there is such a neighborhood which no
geodesic through $x$ intersects in more than one connected component.  Let
$y$ be any point in this neighborhood.  Then there is a unique geodesic
segment inside the neighborhood connecting $x$ to $y$, and no other
geodesic can connect $x$ to $y$ (by the unitrace condition).  Hence $X$
intersects $Y$ in a single point.

Now, $X$ and $Y$ are smooth embedded copies of $\P^{n-1}$ in $G(M)$, and
under smooth deformations the intersection number is constant modulo 2. But
for any $z \in M$, $Y$ can be smoothly deformed into $Z$ (take a smooth
curve in $M$ from $y$ to $z$, and the corresponding skies in $G(M)$ for
each point of this curve).  Hence $X$ intersects $Z$ in a set with an odd
number of elements, and there is always at least one geodesic from $x$ to
$z$.

Finally, since $x$ was arbitrary, any two points in $M$ can be connected by
a geodesic.  ${}_\Box$

\medskip\noindent
The space $(\wM,\wD)$ of Example~\ref{unp} is unitrace and
geodesically connected, but $G(\wM)$ is not Hausdorff.  Thus, for unitrace
spaces $M$, the Hausdorff condition on $G(M)$ is sufficient but not
necessary for geodesic connectedness of $M$.

The preceding theorem is a peer of Proposition 3.5 of \cite{BP2}, in that
there are unitrace spaces with conjugate points ({\it e.g.,} a plane with a
bump), and there are spaces with no conjugate points which are not unitrace
({\it e.g.,} a cone less its vertex).

One can study the geometry of $(M,\D)$ {\em via\/} the topology of $G(M)$.
As an example, let $M$ be an open, convex subset of
$\R^n$ with any linear connection
\D\ whose geodesics are the intersections of the usual straight lines with
$M$.
\begin{proposition}
For any such $M$, $G^+(M) \cong TS^{n-1}$ whence $G(M)$ is
diffeomorphic to an $(n-1)$-plane bundle over $\P^{n-1}$.
\end{proposition}
{\bf Proof:} Choose a point $x_0$ in $M$ and a Euclidean sphere $S$ centered
at $x_0$ of a sufficiently small radius so that it is entirely contained in
$M$. Let \g\ be an oriented geodesic. Euclidean parallel translate \g\ to
$x_0$ and let $x$ be the point of intersection of the translation with $S$ in
the positive direction. Regarding $T_xS$ as embedded in $\R^n$, let $P_x$
denote
the Euclidean orthogonal projection of $M$ onto $T_xS$. Extending \g\ if
necessary, \g\ intersects $P_x$ in a unique point
$y$. Now identify \g\ with $(x,y) \in TS$ and note that $P_x \cong \R^{n-1}$.
${}_\Box$

\medskip\noindent
As a specific case, using the Cayley-Klein model of hyperbolic $n$-space
${\Bbb H}^n$ we have
\begin{corollary}
For any
standard pseudoeuclidean $\R^n$, $G^+({\Bbb H}^n) \cong
G^+({\Bbb R}^n) \cong TS^{n-1}$ whence $G({\Bbb H}^n) \cong G({\Bbb
R}^n)$ as vector bundles over $\P^{n-1}$. ${}_\Box$
\end{corollary}

Geodesic connectedness is a typical example of a geometrically interesting
property. The spaces of the preceding example are geodesically connected, as
are any which are sufficiently like them. In particular,
\begin{theorem}
If $(M,\D)$ has $G^+(M) \cong TS^{n-1}$ with skies
corresponding to sections over $S^{n-1}$, then it is geodesically connected.
\label{ngc}
\end{theorem}
{\bf Proof:} For brevity, let us denote by $B$ the space
$G(M)$ regarded as an $(n-1)$-plane bundle over $\P^{n-1}$.  Now recall
the Stiefel-Whitney classes of a vector bundle; {\em e.g.,} \cite{MS}.  We
shall only be concerned with the top class $w_{n-1}(B)$ which is an element
of $H^{n-1}(\P^{n-1}) = \Z_2$, the integers modulo 2. A standard theorem in
obstruction theory is that if $w_{n-1} \neq 0$, then $B$ has no
nonvanishing sections (if it did, one could split off a line bundle).  We
show $w_{n-1}(B) \neq 0$.

There are two cases.  When $n-1$ is even, $w_{n-1}(B) = w_{n-1}(TP^{n-1})
\neq 0$ because $-I$ is homotopic to $I$ as elements of $GL_{n-1}$ {\em
and\/} as maps of $S^{n-1}$.  When $n-1$ is odd, $w_{n-1}(B) \neq
w_{n-1}(TP^{n-1}) = 0$ because $-I$ is not homotopic to $I$ as elements of
$GL_{n-1}$ {\em and\/} as maps of $S^{n-1}$.  [Indeed, note that in the
first case $\deg(-I) = 0$ on $S^{n-1}$ and in the second $\deg(-I) = -1$.]

Thus $w_{n-1}(B) \neq 0$ in both cases, so $B$ has no nonvanishing
sections, so $M$ is geodesically connected.  ${}_\Box$

\medskip\noindent
In particular, we obtain a new proof of a venerable result.
\begin{corollary}
${\Bbb H}^n$ and any standard pseudoeuclidean $\R^n$ are geodesically
connected. ${}_\Box$
\end{corollary}

More generally, the proof of Theorem \ref{ngc} applies to any $(M,\D)$ such
that $G^+(M)$ is homeomorphic to $TS^{n-1}$ with
the skies
of points in $M$ corresponding to the sections of $TS^{n-1}$ over $S^{n-1}$.
Thus we introduce
\begin{definition}
Let $(E,p,S)$ be a vector bundle.
We say $G(M)$ is a {\em sky bundle\/} over $S$, if there is a
homeomorphism $F: G(M) \rightarrow E$ such that
for each $x \in M$ the
map  $\vp_x = p(F|X): X \rightarrow S$ is a bijection.
\end{definition}
Notice that if $G(M)$ is an $E$-sky bundle over $S$, then the
skies of $M$ correspond to sections of the vector bundle $(E,
p, S)$.  For each $[\g] \in X$ one has that $F[\g]$ is an
element in the fiber over $p(F[\g])$.  Let $h_x: S \rightarrow E$ be the
section corresponding to $x \in M$.  Thus, we have that $h_x(s)
= F( \vp_x^{-1}(s))$.  Let $x$ and $y$ denote distinct points
of $M$.  Then, $h_x - h_y$ vanishes at some $s_0 \in S$ if and only if there
is a geodesic from $x$ to $y$ in $M$.

\begin{lemma}
Let $G(M)$ be a sky bundle $(E, p, S)$.
If $w_{n-1}(E) \not= 0$, then $M$ is geodesically connected.
\end{lemma}
{\bf Proof:} Given any two points $x, y \in M$, let $h_x$ and $h_y$
be the induced sections.  Then $w_{n-1}(E) \not= 0$ implies that
there is some point $s_0 \in S$ with $h_x(s_0) - h_y(s_0) = 0$. ${}_\Box$

\medskip
Let $B^{2n-2} = TS^{n-1}/\sim$ be the
identification described earlier. We regard $TS^{n-1}$ as a subset
of $\R^{2n}$ with coordinates $(p, v)$ for $p$
a (Euclidean) unit vector and $v$ a vector (Euclidean) orthogonal to
$p$. We identify $(p, v) \sim (p, -v)$ to get the bundle $B^{2n-2}$.
 The base space of the vector bundle $B^{2n-2}$ is $\P^{n-1}$ and the
fibers are copies of $\R^{n-1}$.
\begin{corollary}
Let $G(M)$ be a $B^{2n-2}$-sky bundle over $\P^{n-1}$.
Then $M$ is geodesically connected. ${}_\Box$
\end{corollary}

Let $F^+: G^+(M) \rightarrow TS^{n-1}$ be a homeomorphism onto a tubular
neighborhood of the zero section and use the above
notation.  Assume that for each pair of geodesics $\g^+$ and
$\g^-$ corresponding to the same geodesic with opposite
orientations one has $F^+(\g^+) = (p, v)$ if and only if $F^+(\g^-) = (-p,
v)$.  Then there is an induced map $F: G(M) \rightarrow B^{2n-2}$ which
expresses $G(M)$ as a sky bundle over $\P^{n-1}$.  Thus, $M$ is
geodesically connected.

We shall now apply this to Hadamard manifolds. First we generalize a
definition from \cite[p.\,2]{GW}.
\begin{definition}
A manifold with linear connection $(M,\D)$ has a\/ {\em pole @}
if the exponential map\/ $\exp_@: T_@M \rightarrow M$ is a diffeomorphism.
If every point of $M$ is a pole, we call it a\/ {\em pseudohadamard}
manifold.
\end{definition}
Pseudohadamard manifolds share an essential feature of Hadamard
manifolds, as we recall \cite[p.\,1]{GW}.
\begin{definition}
We say that $(M, g)$ is a {\em Hadamard manifold\/} if it is a
simply connected, complete Riemannian (positive-definite) manifold
of nonpositive sectional curvature.
\end{definition}
Greene and Wu \cite{GW} refer to Hadamard manifolds as
Cartan-Hadamard (CH) manifolds.
The universal cover of a complete Riemannian manifold of
nonpositive sectional curvature is a Hadamard manifold.
In a Hadamard manifold, the exponential
map is a diffeomorphism at each point. Thus, two distinct
points of a Hadamard or pseudohadamard manifold are joined by exactly
one geodesic.

\begin{lemma}
A pseudohadamard manifold is nonreturning.
\end{lemma}
{\bf Proof:} As noted above, in a pseudohadamard manifold every
two points are joined by exactly one geodesic segment.
Therefore once a geodesic leaves a convex normal
neighborhood it cannot return. ${}_\Box$
\begin{proposition}
If $(M, \D)$ is a pseudohadamard manifold, then it is pseudoconvex
and $G(M)$ is a manifold.
\end{proposition}
{\bf Proof:} Since $M$ is nonreturning, by Corollary 5.6 of \cite{BP2}
it is sufficient
to show that $M$ is \D-pseudoconvex. First, we note that $M\cong\R^n$ and
$TM\cong M\times\R^n$. We shall regard each $\exp_x:\R^n\rightarrow M$.
Let $K\subseteq M$ be compact and $x\in K$. Then $\exp_x^{-1}(K)
\subseteq\R^n$ is compact. Let $C_x = \{\lambda u ; 0\leq\lambda\leq 1\
\mbox{\ns and}\ u\in\exp_x^{-1}(K)\}$ be the closed cone on
$\exp_x^{-1}(K)$,
so it is also compact. Observe that $H = \bigcup_{x\in K} \exp_x(C_x)$
is the \D-geodesic hull of $K$.

We consider the usual metric on $\R^n$ and choose any
homeomorphism $M\cong\R^n$ to induce a metric on $M$. Let $\delta(x)
= \mbox{\ns diam}(\exp_x(C_x))$. Observe that $x\mapsto C_x$ is a
continuous function from $K$ to the set of compact sets in $\R^n$ ({\em
e.g.,} with respect to the Hausdorff metric).
It follows that $\delta$ is a continuous function $K\rightarrow\R$. Thus it is
bounded
above on $K$, say by $b$. Then $\{y\in M ; d(y,K)\leq b\}$ is a compact
set in $M$ which contains $H$, so $M$ is \D-pseudoconvex. ${}_\Box$

\medskip
For the rest of this paper $(M, g)$ will denote a
Hadamard manifold.  Let a geodesic triangle lying in
$(M, g)$ have sides $c_1, c_2$ and
$c_3$  of lengths $a_1, a_2$,and $a_3$, respectively.
Assume the angle across from the side of length $a_3$ has
measure $\alpha_3$.  For Hadamard manifolds one has the
following cosine inequality:
$$a_3{}^2 \geq a_1{}^2 + a_2{}^2 - 2 a_1 a_2 \cos(
\alpha_3)\,.$$
This is known as the first law of cosines; see
\cite{BGS}, p.\,7.

Let $p$ be a fixed point
of the Hadamard manifold $(M, g)$ and let  $c: \R^1
\rightarrow M$ be a geodesic.  Assume  $q$ is a closest
point to  to $p$ on the image of $c$.  Notice that if $c_1:
[0, a_1] \rightarrow M$ is the unique unit speed geodesic
from $p$ to $q$, then $c_1$ is orthogonal to $c$ at $q$.  If
$c_3: [0, a_3] \rightarrow M$ is a unit speed geodesic
from $p$ to some point $r$ ($r\neq  q$) of $c$, then the above
inequality shows that $a_3 > a_1$.  Hence, $p$ has a unique
closest point $q$ on $c$.  The point $q$ is called the {\em
foot\/} of $p$ on $c$.

Fixing $p$, let $S^{n-1} = \{X \in T_{p}M :  g(X, X) = 1
\}$.  We now define a smooth map  $F: TS^{n-1} \rightarrow
G^{+}(M)$ for $(M, g)$ as follows.  Let $X \in S^{n-1}$
and  $Y \in T_{p}M$ with  $g(X, Y) = 0$.  Then $F(X, Y)$ is
the geodesic with initial point $\exp_{p}(Y)$ and initial
direction the parallel transport of $X$ along the geodesic
$c_{1}(t) = \exp_{p}(tY)$.  The uniqueness of feet for
Hadamard manifolds mentioned above yields that $F$ is a
diffeomorphism.  Thus, we have established the following
result.
\begin{proposition}
If $(M, g)$ is an $n$-dimensional
Hadamard manifold, then $G^{+}(M)$ is diffeomorphic to
$TS^{n-1}$.
\end{proposition}
Once again, we obtain a proof that $\R^n$ and $\H^n$ are geodesically
connected.

\noindent{\bf Acknowledgement} The authors would like to thank
the referee for some helpful comments.

\end{document}